\documentclass[12pt]{amsart}

\usepackage{epsf, amsfonts, amsmath, amsthm, amssymb, latexsym}

\title[Critical {\sf GJMS} operators using Wodzicki's residue]
{A construction of critical {\sf GJMS} operators using Wodzicki's residue}

\author{William J. Ugalde}

\newcommand{\nnn}[1]{\eqref{#1}}

\def\sideremark#1{\ifvmode\leavevmode\fi\vadjust{\vbox to0pt{\vss
\hbox to 0pt{\hskip\hsize\hskip1em
\vbox{\hsize2cm\tiny\raggedright\pretolerance10000
\noindent #1\hfill}\hss}\vbox to8pt{\vfil}\vss}}}

\newtheorem{theorem}{Theorem}[section]
\newtheorem{lemma}[theorem]{Lemma}
\newtheorem{proposition}[theorem]{Proposition}

\newtheorem{definition}[theorem]{Definition}
\newtheorem{remark}[theorem]{Remark}

\newcommand{\bbC}{\mathbb{C}}
\newcommand{\bbR}{\mathbb{R}}
\newcommand{\bbS}{\mathbb{S}}

\newcommand{\cD}{\mathcal{D}}
\newcommand{\cH}{\mathcal{H}}
\newcommand{\cP}{\mathcal{P}}

\renewcommand{\a}{\alpha}
\renewcommand{\b}{\beta}
\renewcommand{\d}{\delta}
\newcommand{\ga}{\gamma}
\newcommand{\Ga}{\Gamma}
\newcommand{\La}{\Lambda}
\newcommand{\Na}{\nabla}
\newcommand{\s}{\sigma}

\newcommand{\Rho}{{\mbox{\sf P}}} % The Schouten tensor

\def\<#1,#2>{\langle{#1},{#2}\rangle} %% my version for inner product
\DeclareMathOperator{\tr}{trace}  %%% uses `amsopn' package
\DeclareMathOperator{\Tr}{Tr} 

\def\set#1{\{\,#1\,\}}           %% set notation
\def\del{\partial}               %% partial derivatives symbol

\def\lcovd#1#2{{#1}_{\,;\,#2}}                      %% for f|_i
\def\ucovd#1#2{{#1}_{\,;}{}^{#2}}                   %% for f|^i

\def\llcovd#1#2#3{{#1}_{\,;\,#2#3}}                     %% for  f|_{ij}
\def\lucovd#1#2#3{{#1}_{\,;\,#2}{}^{#3}}                %% for  f|_i^j
\def\uucovd#1#2#3{{#1}_{\,;}{}^{#2#3}}                  %% for  f|^{ij}

\def\ulucovd#1#2#3#4{{#1}_{\,;}{}^{#2}{}_{#3}{}^{#4}}       %% for  f|^i_j^k
\def\llucovd#1#2#3#4{{#1}_{\,;\,#2#3}{}^{#4}}               %% for  f|_{ij}^k
\def\uuucovd#1#2#3#4{{#1}_{\,;}{}^{#2#3#4}}                 %% for  f|^{ijk}
\def\lllcovd#1#2#3#4{{#1}_{\,;\,#2#3#4}}                    %% for  f|_{ijk}

\def\lulucovd#1#2#3#4#5{#1_{\,;\,#2}{}^{#3}{}_{#4}{}^{#5}}      %% for  f|_i^j_k^l
\def\uulucovd#1#2#3#4#5{#1_{\,;}{}^{#2#3}{}_{#4}{}^{#5}}    %% for  f|^{ij}_k^l

\def\ululucovd#1#2#3#4#5#6{{#1}_{\,;}{}^{#2}{}_{#3}{}^{#4}{}_{#5}{}^{#6}}
\renewcommand{\flat}{\mathrm{flat}}          %% subscript for flat case
\newcommand{\confflat}{\mathrm{conf.\,flat}} %% conformally flat case

\DeclareMathOperator{\Image}{Im}   %% Image
\DeclareMathOperator{\Kernel}{Ker}   %% Kernel
\DeclareMathOperator{\lot}{lot}   %% lower order terms
\DeclareMathOperator{\Rc}{Rc}     %% Ricci tensor
\DeclareMathOperator{\Sc}{Sc}     %% scalar curvature
\DeclareMathOperator{\Wres}{Wres} %% Wodzicki residue
\DeclareMathOperator{\wres}{wres} %% Wodzicki residue density

\begin{document}

%%-----> Abstract:

\begin{abstract}
For an even dimensional, compact, conformal manifold without boundary we construct 
a conformally invariant differential operator of order the dimension of the manifold.
In the conformally flat case, this operator coincides with the critical {\sf GJMS} operator
of Graham-Jenne-Mason-Sparling.  We use the Wodzicki residue of a pseudo-differential operator of 
order~$-2,$ originally defined by A. Connes, acting on middle dimension forms.
\par
\noindent Math.Subj.Clas.: 53A30. Keywords: Wodzicki's residue, GJMS operators. 
\par
\noindent Research supported in part by NSF grant DMS-9983601
\end{abstract}

%\maketitle

%% 53A30 = Conformal differential geometry
%% 35S05 = General theory of PsDO
%% 58J40 = Pseudo-differential and Fourier integral operators on manifolds
%% 58B34 = Noncommutative geometry (à la Connes)
%% 46L87 = Noncommutative differential geometry

\maketitle 

%\tableofcontents

%%-----> INTRODUCTION: ---------------------------------------------------------------------------

\section{Introduction}

In \cite{Con1} Connes uses his quantized calculus to find a conformal invariant.  A central part of the 
explicit computation of this conformal invariant in the 4-dimensional case is the study of a trilinear 
functional on smooth functions over the manifold $M$ given by the relation
$$
\tau(f_0,f,h) = \Wres(f_0[F,f][F,h]),
$$
where $\Wres$ represents the Wodzicki residue, and $F$ is a pseudo-differential operator of order 0 
acting on 2-forms over $M.$ $F$ is given as $2\cD-1$ with $\cD$ the orthogonal projection on the image 
of $d$ in $\cH=\La^2(T^*M)\ominus H^2$ with $H^2$ the space of harmonic forms.  Using the general 
formula for the total symbol of the product of two pseudo-differential operators, Connes computed 
a natural bilinear differential functional  of order $4,$ $B_4$ acting on $C^\infty(M).$  This 
bilinear differential functional 
$$
B_4(f,h) = 4\Delta(\<df,dh>) -2\Delta f \Delta h + 4\<\Na df, \Na dh> +8\<df,dh>J,
$$
is symmetric: $B_4(f,h) = B_4(h,f),$ conformally invariant:
$\widehat{B_4}(f,h) = e^{-4\eta}B_4(f,h),$ for $\widehat{g} = e^{2\eta}g,$
and uniquely determined, for every $f_i$ in $C^\infty(M),$ by the relation:
$$
\tau(f_0,f_1,f_2) = \int_M f_0 B_4(f_1,f_2)\,d  x.
$$
On an $n$ dimensional manifold, we use $J$ to represent the normalized scalar curvature 
$2(n-1)J=\Sc.$

The {\sf GJMS} operators \cite{GJMS} are invariant operators on conformal densities
$$
P_{2k}:\mathcal{E}[-n/2+k] \to \mathcal{E}[-n/2-k]
$$ 
with principal parts $\Delta^k$, unless the dimension is even and $2k > n.$  If $n$ is 
even, the $n$-th  order operator 
$$
P_n:\mathcal{E}[0] \to \mathcal{E}[-n]
$$
is called \emph{critical} {\sf GJMS} operator.  As noted in \cite{GJMS}, 
$\mathcal{E}[0]=C^\infty(M)$ and $\mathcal{E}[-n]$ is the bundle of volume densities on $M.$ 
In the 4-dimensional case, Connes has shown that 
\begin{align*} 
P_4(f) &= \Delta^2 f +2 \Delta(f)\,J + 4\<\Rho,\Na df> + 2\<dJ,df>    \\
       &= \d (d\d +(n-2) J - 4\Rho)d\, f,
\end{align*} 
the Paneitz operator (\emph{critical} {\sf GJMS} for $n=4$), can be derived from $B_4$ by the relation
$$
\int_M B_4(f,h)\,d  x = \frac{1}{2} \int_M f P_4(h)\,d x.
$$

The GJMS operators $P_{2k},$ of Graham-Jenne-Mason-Sparling \cite{GJMS} by construction
have the following properties.
\begin{itemize}
\item[o.] 
$P_{2k}$ exists for all $k$ if $n$ is odd, and if $n$ is even exists for $1\leq k \leq n/2.$
\item[i.] 
$P_{2k}$ is formally self-adjoint.
\item[ii.] 
$P_{2k}$ is conformally invariant in the sense that
$$
\widehat P_{2k}=e^{-(n/2+k)\eta} P_{2k} e^{(n/2-k)\eta}
$$
for conformally related metrics $\widehat g = e^{2\eta}g.$
\item[iii.] 
$P_{2k}$ has a polynomial expression in $\Na$ and $R$ in which all coefficients are rational 
in the dimension $n$.
\item[iv.] 
$P_{2k}=\Delta^k + \lot$.  (Here and below, $\lot = $``lower order terms''.)
\item[v.] 
$P_{2k}$ has the form
$$
\d S_{2k} d + \bigl(\dfrac{n}2-k\bigr) Q_{2k}, 
$$
where $Q_{2k}$ is a local scalar invariant, and $S_{2k}$ is an operator on 1-forms of the form
$$
(d\delta)^{k-1}+\lot \mbox{~or~}\Delta^{k-1}+\lot.
$$
\end{itemize}
In this last expression, $d$ and $\d$ are the usual de Rham operators and $\Delta$ is the form Laplacian 
$d \delta + \delta d$.  (o), (ii), and (iv) are proved in \cite{GJMS}.  The fact that (o) is 
exhaustive is proved in \cite{GoHi}. (i), (iii), and (v) are proved in \cite{Bra5}.

The original work \cite{GJMS} uses the ``ambient metric construction'' of \cite{FeGra-amc} to prove 
their existence.  Since their appearance, there has been a considerable interest in constructive 
ways to obtain formally selfadjoint, conformally invariant powers of the Laplacian on functions with 
properties as stated.  Examples can be found in \cite{GraZw} with scattering theory, in \cite{FeGra} 
with Poincar\'e metrics, and in \cite{GoPe} with tractor calculus.  

In \cite{mio1} (see also \cite{mio2}) we have proved the following result:
\begin{theorem}
\label{theoremonOmeganinmio1}
Let $M$ be a compact conformal manifold without boundary. Let $S$ be a pseudo-differential operator of 
order 0 acting on sections of a vector bundle over $M$ such that $S^2 f_1 = f_1 S^2 $ and the 
pseudo-differential operator $P=[S,f_1][S,f_2]$ is conformal invariant for any $f_i \in C^\infty(M).$  
There exists a unique, symmetric, and conformally invariant, bilinear, diffe\-ren\-tial functional 
$B_{n,S}$ of order $n$ such that
$$
\Wres(f_0[S,f_1][S,f_2]) = \int_M f_0 B_{n,S}(f_1,f_2)\, d x
$$
for all $f_i \in C^\infty(M).$
\end{theorem}

This result, aiming to extend part of the work of Connes~\cite{Con1}, is based on the study of the 
formula for the total symbol $\sigma(P_1 P_2)$ of
the product of two pseudo-differential operators $P_1$ and $P_2$ given by
$$
\sum \frac{1}{\a!} \del^\a_\xi(\sigma(P_1)) \,D^\a_x(\sigma(P_2))
$$
where $\a=(\a_1,\dots,\a_n)$ is a multi-index, $\a! = \a_1!\cdots \a_n!$ and 
$D^\a_x = (-i)^{|\a|} \del^\a_x.$ This expression for $\s(P_1P_2)$ is highly asymmetric on $\del_\xi$ 
and $\del_x.$ Even though, some control can be achieved when $P_1$ say, is given as a ``multiplication'' 
operator.

If $F$ is the operator $2\cD-1$ where $\cD$ is the orthogonal projection on  the image of $d$ inside 
the Hilbert space of square integrable forms of middle dimension without the  harmonic ones, $\cH = 
L^2(M,\Lambda^{n/2}_{\bbC} T^*M)\ominus H^{n/2}$ as defined in \cite{Con1} then, this differential 
functional $B_n=B_{n,F}$ provides a way of constructing operators like the critical {\sf GJMS} 
operators on compact conformal manifolds of even dimension, out of the Wodzicki residue of a 
commutator operator acting on middle dimension forms.  We say here `like' because up to date, we do 
not know the relation in between the operators we construct and the ``ambient metric construction'' 
of \cite{FeGra-amc}.  The following is the main result of this paper:
\begin{theorem}
\label{theoremonP_n}
Let $M$ be an even dimensional, compact, conformal manifold without boundary and let $(\cH,F)$ be 
the Fredholm module associated to $M$ by Connes \cite{Con1}. Let $P_n$ be the differential 
operator given by the relation
\begin{equation}
\label{omegatoP}
\int_M B_n(f,h) \,d x = \int_M f P_n(h) \,d x
\end{equation}
for all $f,h \in  C^\infty(M).$  Then,
\begin{itemize}
\item[i.] 
$P_n$ is formally selfadjoint.
\item[ii.] 
$P_n$ is conformally invariant in the sense $\widehat{P_n}(h) = e^{-n\eta} P_n(h),$ if
$\widehat g = e^{2 \eta}g.$
\item[iii.] 
$P_n$ is expressible universally as polynomial in the ingredients $\Na$ and $R$ with 
coefficients rational in $n.$
\item[iv.] 
$P_n(h) = c_n \Delta^{n/2}(h) + \lot$ with $c_n$ a universal constant and $\lot$
meaning ``lower order terms''.
\item[v.] 
$P_n$ has the form $\delta S_n d$ where $S_n$ is an operator on 1-forms given as
a constant multiple of $\Delta^{n/2-1} + \lot,$ or $(d\delta)^{n/2-1} + \lot.$
\item[vi.] 
$P_n$ and $B_n$ are related by:
$$ 
P_n(fh) - fP_n(h) - hP_n(f) = -2B_n(f,h).
$$
\end{itemize}
\end{theorem}

The structure of this paper is as follows. In Section~\ref{notations} we provide the setting, 
definitions, and conventions we will use throughout the rest of the paper plus a small explanation of 
the origin of $\Wres(f_0[F,f][F,h]).$   In  Section~\ref{bilineardfiffunc} we review 
Theorem~\ref{theoremonOmeganinmio1} to make this paper 
as self contained as possible.  In Section~\ref{thesymbolofF}, by a recursive computation of the 
symbol expansion of $F,$ we show that $B_n(f,h)$ has a universal expression as a polynomial on 
$\Na,$ $R,$ $df,$ and $dh.$  In Section~\ref{PnandWres} we prove Theorem~\ref{theoremonP_n}.  In 
\cite{Dani} we will present the computations conducting to the explicit expression in 
the 6-dimensional case together with other calculations related to subleading symbols of 
pseudo-differential operators.

Special thanks to Thomas Branson for support and guidance, to the 
Erwin Schr\"odinger International Institute for Mathematical Physics for a pleasant environment 
where some of the last stages of this work took place, and also to Joseph V\'arilly for suggestion 
on how to improve an earlier version. 

%%-----> NOTATION SECTION: ---------------------------------------------------------------------

\section{Notations and conventions}
\label{notations}

In this paper, $R$ represents the Riemann curvature tensor, $\Rc_{ij} = R^{k}{}_{ikj}$ represents 
the Ricci tensor, and $\Sc =\Rc^{i}{}_{i}$ the scalar curvature. $\Rho = (\Rc - Jg)/(n - 2)$ is the 
Schouten tensor.  The relation between the Weyl tensor and the Riemann tensor is given by
$$
W^i{}_{jkl} = 
R^i{}_{jkl} + \Rho_{jk}{\d_K}^i{}_l - \Rho_{jl}{\d_K}^i{}_k + \Rho^i{}_l g_{jk} - \Rho^i{}_k g_{jl}
$$
where $\d_K$ represents the Kronecker delta tensor.  If needed, we will ``raise'' and ``lower'' indices
without explicit mention, for example, $g_{mi}R^i{}_{jkl}=R_{mjkl}$.

A differential operator acting on $C^\infty(M)$ is said to be natural if it can be written as a 
universal polynomial expression in the metric $g,$ its inverse $g^{-1},$ the connection $\Na,$ and the 
curvature $R;$ using tensor products and contractions. The coefficients of such natural operators 
are called natural tensors.  If $M$ has a metric $g$ and if $D$ is a natural differential operator 
of order $d$, then $D$ is said to be conformally invariant with weight $w$ if after a conformal 
change of the metric $\widehat{g} = e^{2\eta} g$, $D$ transforms as 
$\widehat{D} = e^{-(w+d) \eta} D e^{w\eta}.$  In the particular case in which $w=0$ we say that
$D$ is conformally invariant and the previous reduces to $\widehat{D} = e^{-d\eta} D.$
For a natural bilinear differential operator $B$ of order $n,$ conformally invariance with bi-weight
$(w_1,w_2)$ means
\begin{equation}
\label{confinv1}
\widehat{B}(f,h) = e^{-(w_1+w_2+n)\eta} B(e^{w_1\eta}f,e^{w_2\eta}h)
\end{equation}
for every $f,h \in C^\infty(M).$  In this work will consider only the case $w_1=w_2=0.$
In this particular situation, because $\widehat{d x} = e^{n\eta} d x$ we have
\begin{equation}
\label{confinv2}
\widehat B(f,h)\,\widehat{d x} =B(f,h)\,d x.
\end{equation}
Some times, when the context allows us to do so, we will abuse of the language and say
that $B$ is conformally invariant making reference to \nnn{confinv2} instead of \nnn{confinv1}.

When working with a pseudo-differential operator $P$ of order $k,$ its total symbol 
(in some given local coordinates) will be represented as
a sum of $r \times r$ matrices of the form $\s(P) \sim \s_k^P + \s_{k-1}^P + \s_{k-2}^P+\cdots,$
where $r$ is the rank of the vector bundle $E$ on which $P$ is acting.
It is important to note that the different $\s_j^P$ are defined only in local charts and
are not diffeomorphism invariant \cite{Ku}.  However, Wodzicki \cite{Wo} has shown that the term
$\s_{-n}^P$ enjoys a very special significance.   For a pseudo-differential operator $P,$ acting on
sections of a bundle $E$ over a manifold $M,$ there is a $1$-density on $M$ expressed in local
coordinates by
\begin{equation}
\label{wresP}
\wres(P) = \int_{||\xi||=1}\tr\bigl(\s_{-n}^P(x,\xi)\bigr)\,d \xi \,d x,
\end{equation}
where $\s^P_{-n}(x,\xi)$ is the component of order $-n$ in the total symbol of $P,$ $||\xi||=1$
means the Euclidean norm of the coordinate vector $(\xi_1,\cdots,\xi_n)$ in $\bbR^n,$ and
$d \xi$ is the normalized volume on $\set{||\xi||=1}.$ This \emph{Wodzicki residue density} is
independent of the local representation. An elementary proof of this matter can be found in 
\cite{Fedosov}.  The Wodzicki residue, $\Wres(P)=\int \wres(P)$ is then independent of the choice of 
the local coordinates on $M$, the local basis of $E$, and defines a trace (see \cite{Wo}).

For a compact, oriented manifold $M$ of even dimension $n=2l,$ endowed with a conformal structure,
there is a canonically associated Fredholm module $(\cH,F)$ \cite{Con1}.  
$\cH$ is the Hilbert space of square integrable forms of middle 
dimension with an extra copy of the harmonic forms 
$\cH =L^2(M,\Lambda^{n/2}_{\bbC} T^*M)\oplus H^{n/2}$.  Functions on $M$ act as multiplication 
operators on $\cH.$ $F$ is a pseudo-differential operator of order 0 acting in $\cH.$ On $\cH 
\ominus H^{n/2}\ominus H^{n/2} $ it is obtained from 
the orthogonal projection $\cD$ on the image of $d,$ by the relation $F=2\cD-1.$  On 
$H^{n/2} \oplus H^{n/2}$ it only interchanges the two components.
We are only concerned with the non-harmonic part of $\cH.$  From the Hodge decomposition theorem  it is 
easy to see that in terms of a Riemannian metric compatible with the conformal structure of $M$, $F$ 
restricted to $\cH\ominus H^{n/2}\ominus H^{n/2}$ can be written as:
$$
F = \frac{d \d-\d d}{d \delta + \delta d}.
$$

%% Why $\Wres(f_0[F,f][F,h])$?: ------------------------------------------------------------------

\subsection{Why $\Wres(f_0[F,f][F,h])$?}

As pointed out in \cite{OPS}, the quantum geometry of strings is concerned, among other ideas, with 
determinants of Laplacians associated to a surface with varying metrics. For a compact surface $M$ 
without boundary and metric $g,$ it is possible to associate to its Laplacian $\Delta_g$ a determinant
$$
\det{}' \Delta_g = \prod_{\lambda_j \not= 0} \lambda_j
$$
where $0 \leq \lambda_0 < \lambda_1 \leq \lambda_2 \leq \cdots$ are the eigenvalues of $\Delta_g$.

To make sense out of this formal infinite product, some regularization procedure is needed. The key
to study this expression as a function of $g$ is the \emph{Polyakov action} \cite{Pol} which gives
a formula for the variation of $\log \det' \Delta_g$ under a conformal change of $g.$
For a Riemann surface $M$, a map $f = (f^i)$ from $M$ to $\bbR^2$ and metric $g_{ij}(x)$ on $M$, the
2-dimensional Polyakov action is given by
$$
I(f) = \frac{1}{2\pi} \int_M g_{ij} \,df^i \wedge \star d f^j.
$$

By considering instead of $d f$ its quantized version $[F,f],$ Connes \cite{Con1} quantized the Polyakov 
action as a Dixmier trace:
\begin{align*}
\frac{1}{2 \pi} \int_M g_{ij} d f^i \wedge \star d f^j = -\frac{1}{2}
\Tr_\omega\bigl(g_{ij}[F,f^i][F,f^j]\bigr).
\end{align*}
Because Wodzicki's residue extends uniquely the Dixmier trace as a trace on the
algebra of PsDOs, this quantized Polyakov action has sense in the general even dimensional case.

Connes' trace theorem \cite{ConAction} states that the Dixmier trace and the Wodzicki residue of an 
elliptic PsDO of order $-n$ in an $n$-dimensional manifold $M$ are proportional by a factor of 
$n(2\pi)^n.$  In the 2-dimensional case the factor is $8\pi^2$ and so the quantized Polyakov action 
can be written as,
\begin{align*}
-16\pi^2 I=\Wres \bigl(g_{ij} [F,f^i][F,f^j]\bigr).
\end{align*}

%%-----> A bilinear differential functional: --------------------------------------------------

\section{A bilinear differential functional}
\label{bilineardfiffunc}

For a vector bundle $E$ of rank $r$ over a compact manifold without boundary $M$ of dimension $n,$
let $S$ be a pseudo-differential operator of order $k$ acting on sections of $E.$ Consider $P$ as 
the pseudo-differential operator given by the product $P= f_0[S,f][S,h]$ with each $f_0,f,h \in 
C^\infty(M)$.
$P$ is acting on the same vector bundle as $S$, where smooth functions on $M$ act as multiplication 
operators. Each commutator $[S,f]$ defines a pseudo-differential operator of order $k-1,$ thus $P$ has
order $2k-2$. Assuming $2k-2 \geq -n,$ in a given system of local coordinates the total symbol of $P$, 
up to order $-n$, is represented as a sum of $r \times r$ matrices of the form 
$\s_{2k-2}^P+\s_{2k-3}^P+\cdots+\s_{-n}^P.$  We aim to study
$$
\Wres(P) = \int_M \int_{||\xi||=1}\tr\bigl(\s^P_{-n}(x,\xi)\bigr)\,d \xi\,d x.
$$
Unless otherwise stated, we assume a given system of local coordinates.

For the operator multiplication by $f$ we have $\s(f) = f I$ with $I$ the identity operator on the 
sections on which it is acting. In this way $\s(f P_2) = f\s(P_2)$ and in particular 
$\s_{-n}(f P_2) = f\s_{-n}(P_2).$ 
As a consequence we obtain the relation
$$
\wres(f_0 P_2) = f_0 \wres(P_2).
$$

The proof of the following two results is similar to the proof for the case $S$ of 
order 0 in Lemma~2.2 and Lemma~2.3 in~\cite{mio1} so we omit them in here.

\begin{lemma}
$[S,f]$ is a pseudo-differential operator of order $k-1$ with total symbol
$\s([S,f]) \sim \sum_{j \geq 1} \s_{k-j}([S,f])$ where
\begin{equation}
\label{sym-nSf} 
\s_{k-j}([S,f]) =\sum_{|\b|=1}^j \frac{D^\b_x(f)}{\b!} \del_\xi^\b \bigl(\s^S_{k-(j-|\b|)}\bigr).
\end{equation}
\end{lemma}

\begin{lemma}
\label{lemmas_nSS} For $2k+n \geq 2,$ with the sum taken over $|\a'|+|\a''|+|\b|+|\d|+i+j = n+2k$, 
$|\b|\geq 1$, and $|\d|\geq 1$,
$$
\s_{-n}([S,f][S,h]) 
=\sum \frac{D^\b_x(f) D_x^{\a'+\d}(h)}{\a'!\a''!\b!\d!}
   \del_\xi^{\a'+\a''+\b}(\s^S_{k-i}) \del_\xi^\d(D_x^{\a''}(\s^S_{k-j})).
$$
\end{lemma}

\begin{definition}
\label{defB}
For every $f,h \in C^\infty(M)$ we define $B_{n,S}(f,h)$ as given by the relation
$$
B_{n,S}(f,h)\,d x := \wres([S,f][S,h]).
$$
\end{definition}

\noindent Because of Lemma~\ref{lemmas_nSS}, $B_{n,S}(f,h)$ is explicitly given by:
$$
\int_{||\xi||=1} 
\sum \frac{D^\b_x(f) D^{\a''+\d}_x(h)}{\a'!\a''!\b!\d!}
\tr\biggl(
\del^{\a'+\a''+\b}_\xi(\s^S_{k-i}) \del^{\d}_\xi(D^{\a'}_x(\s^S_{k-j})) \biggr)\,d \xi
$$
with the sum taken over $|\a'|+|\a''|+|\b|+|\d|+i+j = n+2k$,
$|\b|\geq 1$, and $|\d|\geq 1.$

The arbitrariness of $f_0$ in the previous construction, implies that $B_{n,S}$ is
uniquely determined by its relation with the Wodzicki residue of the operator $f_0[S,f][S,h]$
as stated in the following theorem.

\begin{theorem}
\label{theoremonomegan1} 
There is a unique bilinear diffe\-ren\-tial functional $B_{n,S}$ of order $n$ such that
$$
\Wres(f_0[S,f][S,h]) = \int_M f_0 B_{n,S}(f,h) \,d x
$$
for all $f_0,$ $f,$ $h$ in $C^\infty(M).$  Furthermore, as the left hand side, the right hand side 
defines a Hochschild $2$-cocycle over the algebra of smooth functions on $M.$
\end{theorem}

\begin{proof}
Linearity is evident and uniqueness follows from the arbitrariness of $f_0.$  Let $b$ denote the 
Hochschild coboundary operator on $C^\infty(M)$ and consider 
$$
\varphi(f_0,f_1,\cdots,f_r) = \Wres(f_0[S,f_1]\cdots[S,f_r]).
$$
From the relation $[S,fh] = [S,f]h + f[S,h]$ we have:
\begin{align*}
&(b\varphi)(f_0,f_1,\cdots,f_r) = \Wres(f_0 f_1 [S,f_2] \cdots [S,f_r]) \\
&\quad + \sum_{j=1}^{r-1}(-1)^j\Wres(f_0 [S,f_1] \cdots [S,f_j f_{j+1}] \cdots [S,f_r]) \\
&\quad      +(-1)^r \Wres(f_rf_0[S,f_1] \cdots [S,f_{r-1}]) \\
&=\Wres(f_0 f_1 [S,f_2] \cdots [S,f_r]) \\
&\quad + \sum_{j=2}^{r}(-1)^{j-1}\Wres(f_0 [S,f_1] \cdots [S,f_{j-1}] f_j [S,f_{j+1}] \cdots [S,f_r]) 
\\
&\quad + \sum_{j=1}^{r-1}(-1)^{j}\Wres(f_0 [S,f_1] \cdots [S,f_{j-1}] f_j [S,f_{j+1}] \cdots [S,f_r]) 
\\
&\quad + (-1)^r \Wres(f_rf_0[S,f_1] \cdots [S,f_{r-1}]) \\
=& \Wres(f_0 f_1 [S,f_2] \cdots [S,f_r]) 
+ (-1)^{r-1}\Wres(f_0 [S,f_1] \cdots [S,f_{r-1}] f_r) \\
& - \Wres(f_0 f_1 [S,f_2] \cdots [S,f_r]) 
+(-1)^r \Wres(f_rf_0[S,f_1] \cdots [S,f_{r-1}]) = 0
\end{align*}
because of the trace property of $\Wres$.  The result follows as the particular case $r=3.$
\end{proof}

So far, by taking $f_0 = 1$, by uniqueness, and by the trace property of $\Wres$, we conclude
$$ 
\int_M B_{n,S}(f,h) \, d x = \Wres([S,f][S,h]) =  \Wres([S,h][S,f]) = \int_M B_{n,S}(h,f)\, d x.
$$
From Definition~\ref{defB}, $B_{n,S}(f,h)\,d x = \wres([S,f][S,h])$ but in general $\wres$ is not
a trace, hence asserting that $B_{n,S}(f,h)$ is symmetric is asserting that $$
\wres([S,f][S,h])=\wres([S,h][S,f]). $$ To conclude the symmetry of $B_{n,S}(f,h)$ on $f$ and
$h$, it is necessary to request more properties on the operator $S$.  It is enough
to have the property $S^2f=fS^2,$ for every $f\in C^\infty(M).$ The symmetry of $B_{n,S}$ follows from 
the trace properties of the Wodzicki residue and the commutativity of the algebra $C^\infty(M).$  For a 
proof of the following result see \cite{mio1}.

\begin{theorem}
\label{symmetryofOmega} If $S^2 f = f S^2$ for every $f \in C^\infty(M)$ then the differential
functional $B_{n,S}$ in Theorem~\ref{theoremonomegan1} is symmetric in $f$ and $h$.
\end{theorem}

In the most of the present work we will use the particular case $S^2= I.$

%%-----> $B_{n,S}$ on conformal manifolds:------------------------------------------------------

\subsection{The functional on conformal manifolds}
\label{omegaandfredmod}

If we endowed the manifold $M$ of a conformal structure, and ask the operator $P=[S,f_1][S,f_2]$ to 
be independent of the metric in the conformal class, then we can say a little more about this differential 
functional $B_{n,S}$.

\begin{theorem}
\label{theoremonOmeganS} 
Let $M$ be a compact conformal manifold without boundary. Let $S$ be a pseudo-differential operator 
acting on sections of a vector bundle over $M$ such that $S^2 f_1 = f_1 S^2 $ and the 
pseudo-differential operator $P=[S,f_1][S,f_2]$ is conformal invariant for any $f_i \in C^\infty(M).$  
There exists a unique, symmetric, and conformally invariant bilinear diffe\-ren\-tial functional 
$B_{n,S}$ of order $n$ such that
$$
\Wres(f_0[S,f_1][S,f_2]) = \int_M f_0 B_{n,S}(f_1,f_2)\, d x
$$
for all $f_i \in C^\infty(M).$ Furthermore, $\int_M f_0 B_{n,S}(f_1,f_2)\,d x$ defines a
Hochschild 2-cocycle on the algebra of smooth functions on $M.$
\end{theorem}

\begin{proof}
Uniqueness follows from \nnn{theoremonomegan1}.  Symmetry follows from Theorem~\ref{symmetryofOmega} and
its conformal invariance, $\widehat{B_{n,S}} = e^{-n\eta} B_{n,S},$ follows from its 
construction. Indeed, the only possible metric dependence is given by the operator
$P.$ Roughly speaking, the cosphere bundle $\bbS^*M,$ as a submanifold of $T^*M,$ depends on a choice of 
metric in $M.$ But since $\s_{-n}(P)$ is homogeneous of degree $-n$ in $\xi,$ $\{||\xi||=1\}$ can be 
replaced by any sphere with respect to a chosen Riemannian metric on $M,$ $\{||\xi||_g = 1\}.$ 
The formula for the change of variable shows that
$$
\int_{||\xi||=1} \tr(\s_{-n}(P)) d  \xi
$$
does not change within a given conformal class except trough a metric dependence of $P$
(this argument is taken from the proof of Theorem~2.3~\cite{Cordelia}).
The Hochschild cocycle property follows from Theorem~\ref{theoremonomegan1}.
\end{proof}

We restrict ourselves to the particular case of an even dimensional, compact, oriented, conformal
manifold without boundary $M$, and $(E,S)$ given by the canonical Fredholm module $(\cH,F)$ associated
to $M$ by Connes \cite{Con1}, the pseudo-differential operator $F$ of order 0 is given by
$F = (d \d - \d d)(d \d + \d d)^{-1}$ acting on the space
$L^2(M,\La^{l}_\bbC T^*M)\ominus H^{n/2},$ with $H^{n/2}$ the finite dimensional space of middle 
dimension harmonic forms.  By definition $F$ is selfadjoint and such that $F^2=1.$ We relax the notation 
by denoting $B_n = B_{n,F}$ in this particular situation. $P$ is the pseudo-differential operator of 
order $-2$ given by the product $P= f_0[F,f_1][F,f_2]$ with each $f_i \in C^\infty(M)$. $P$ is acting on 
middle dimension forms.  In this situation, Theorem~\ref{theoremonOmeganS} is stated as follows: 
\begin{theorem}
\label{theoremonOmegan}
If $M$ is an even dimensional, compact, conformal manifold without boundary and $(\cH,F)$ is the Fredholm 
module associated to $M$ by Connes \cite{Con1} then, there is a unique, symmetric, and conformally 
invariant bilinear diffe\-ren\-tial functional $B_n$ of order $n$ such that
$$
\Wres(f_0[F,f_1][F,f_2]) = \int_M f_0 B_n(f_1,f_2) \,d x
$$
for all $f_i \in C^\infty(M).$  Furthermore, $\int_M f_0 B_n(f_1,f_2)\, d x$ defines a
Hochschild 2-cocycle on the algebra of smooth functions on $M.$
\end{theorem}

\begin{proof} We must only prove its conformal invariance.
In Lemma~2.9~\cite{mio1} we show that the Hodge star operator restricted to middle dimension forms is 
conformally invariant in fact, acting on $k$-forms we have $\widehat \star = e^{(2k-n)\eta}\star$ for 
$\widehat g= e^{2\eta}g.$ 
The space $\La^{n/2} T^* M$ of middle dimension forms has an inner product
$$
\<\xi_1,\xi_2> = \int_M \overline{\xi_1} \wedge \star \xi_2
$$
which is unchanged under a conformal change of the metric, that is, its Hilbert space completion
$L^2(M,\La^{n/2}T^* M)$ depends only on the conformal class of the metric. Furthermore, the Hodge 
decomposition:
\begin{align*}
\La^{n/2} T^* M &= \Delta(\La^{n/2} T^* M)\oplus H^{n/2}
\\
&=d\delta (\La^{n/2} T^* M)\oplus \delta d(\La^{n/2} T^* M)\oplus H^{n/2},
\end{align*}
is preserved under conformal change of the metric.
Indeed, the space of middle dimension harmonic forms given as $H^{n/2}=\Kernel(d)\cap\Kernel(d\star)$ 
is conformally 
invariant.
Now if $\omega$ is a middle dimension form orthogonal to the space of harmonic forms we have 
$\omega=(d\widehat \delta + \widehat \delta d)\omega_1 = (d\delta + \delta d)\omega_2$ 
with $\omega_1, \omega_2$ middle dimension forms.  It is not difficult to check that for a k-form $\xi$ 
we have:
$$
\widehat \delta \xi = e^{(2(k-1)-n)\eta}\delta e^{(-2k+n)\eta}\xi,
$$
so we must have
$$
d \delta \omega_2 + \delta d \omega_2 =
\omega=(d\widehat \delta + \widehat \delta d)\omega_1 
= d(e^{-2\eta}\delta \omega_1) + \delta(e^{-2\eta}d\omega_1) 
$$
which implies
$$
\omega_0 :=
\underbrace{d(e^{-2\eta}\delta \omega_1 - \delta \omega_2)}_{\in\,(\Image \delta d)^\perp} 
= \underbrace{\delta(d\omega_2 - e^{-2\eta} d\omega_1)}_{\in\,(\Image d\delta)^\perp}
\in H^{n/2}. 
$$
But also $\omega_0$ is orthogonal to any harmonic form so $\omega_0=0$ and hence
$$
d \widehat \delta \omega_1 = d \delta \omega_2
\qquad
\hbox{and}
\qquad
\widehat \delta d \omega_1 = \delta d \omega_2.
$$
We are concerned with $F$ acting on $\Delta(\La^{n/2})=\widehat{\Delta}(\La^{n/2}).$
If
$$
\omega = d\widehat\delta \omega_1 + \widehat\delta d \omega_1
       = d\delta \omega_2 + \delta d \omega_2 ,
$$ 
then 
$$
\widehat{F} \omega = d\widehat\delta \omega_1 - \widehat\delta d \omega_1
=d\delta \omega_2 - \delta d \omega_2 = F \omega
$$
thus $\widehat{F} = F.$

Last, $\Wres([F,f][F,h])$ does not depend on the choice of the metric in the conformal class 
and the result follows.
\end{proof}

%%-----> The functional in the Flat Case:

\subsection{The functional in the flat case}
\label{particularocurrence}

In the particular case of a flat metric, we have $\s_{-k}^F = 0$ for all $k \geq 1$ and so the
symbol of $F$ coincides with its principal symbol given by
$$
\s^F_0(x,\xi) = (\varepsilon_\xi \iota_\xi - \iota_\xi \varepsilon_\xi)||\xi||^{-2}
$$
for all $(x,\xi) \in T^*M$, where $\varepsilon_\xi$ denotes exterior multiplication and $\iota_\xi$ 
denotes interior multiplication.  Using this result, it is possible to give a formula for $B_n$ in the 
flat case using the Taylor expansion of the function
$$
\psi(\xi,\eta)=\tr\bigl(\s_0^F(\xi) \s_0^F(\eta)\bigr).
$$
Here  $\xi,\eta \in T^*_x M \smallsetminus \{0\}$ and  $\s_0^F(\xi) \s_0^F(\eta)$ is acting on
middle-forms.  

\begin{proposition}
\label{theoremtraces0s0}
With $\s^F_0(\xi) \s^F_0(\eta)$ acting on $m$-forms we have:
$$
\psi(\xi,\eta) = \tr(\s^F_0(\xi) \s^F_0(\eta))
= a_{n,m} \frac{\<{\xi},{\eta}>^2}{||\xi||^2||\eta||^2} + b_{n,m},
$$
where
$$
b_{n,m} = \binom{n-2}{m-2} + \binom{n-2}{m} - 2 \binom{n-2}{m-1} = \binom{n}{m} - a_{n,m}.
$$
\end{proposition}
\noindent For the proof see Theorem~4.3~\cite{mio1}.  For the case $n = 4$ and $m = 2$ we obtain
$$
\psi(\xi,\eta) = 8 \<\xi,\eta>^2(||\xi||||\eta||)^{-2} - 2.
$$
Here there is a discrepancy with Connes, he has 
$4 \<\xi,\eta>^2(||\xi||||\eta||)^{-2} +~\mathrm{constant}$. For $n = 6$, $m = 3$,
$$
\psi(\xi,\eta) = 24 \<\xi,\eta>^2(||\xi||||\eta||)^{-2} - 4.
$$

In the flat case, with $S=F,$ \nnn{sym-nSf} reduces to:
$$
\s_{-r}([F,f]) = \sum_{|\b|=r} \frac{D_x^\b f}{\b!} \del_\xi^\b(\s_0^F).
$$
Using this information we deduce from Lemma~\ref{lemmas_nSS} with $\a' = 0,$ $\a'' = \a,$
and by the definition of $B_{n,S}:$
$$
B_{n\,\flat}(f,h) = \int_{||\xi||=1} 
\sum \frac{D^\b_x f D^{\a+\d}_x h}{\a!\b!\d!}
\tr\biggl(
 \del^{\a+\b}_\xi(\s_0^F)\del^\d_\xi(\s_0^F) \biggr)
 \,d \xi
$$
with the sum taken over $|\a|+|\b|+|\d| = n$, $1\leq |\b|$, $1\leq |\d|.$

\begin{definition}
We denote by $T'_n\psi(\xi,\eta,u,v)$ the term of order $n$ in the
Taylor expansion of $\psi(\xi,\eta)$ minus the terms with only powers
of $u$ or only powers of $v$. That is to say,
\begin{equation}
\label{T'}
T'_n \psi(\xi,\eta,u,v) :=
\sum \frac{u^\b}{\b!} \frac{v^\d}{\d!}
     \tr\bigl(\del^\b_\xi(\s_0^F(\xi)) \del^\d_\eta(\s_0^F(\eta))\bigr),
\end{equation}
with the sum taken over $|\b|+|\d|=n,$ $|\b|\geq 1,$ and $|\d|\geq 1.$ 
\end{definition}

There is an explicit expression for $B_n$ in the flat case in terms of the Taylor expansion of 
$\psi(\xi,\eta):$

\begin{theorem}
\label{lemmaomeganflat}
$$
B_{n\,\flat}(f,h)
= \sum A_{a,b} (D_x^a f)(D_x^b h),
$$
where
$$
\sum A_{a,b} u^a v^b =
\int_{||\xi||=1} \bigl(T'_n\psi(\xi,\xi,u+v,v) -
                     T'_n\psi(\xi,\xi,v,v) \bigr) \,d \xi
$$
with $T'_n\psi(\xi,\eta,u,v)$ is given by \nnn{T'}.
\end{theorem}

For details and a proof of the previous theorem see Section~4~\cite{mio1}.  In the 4 dimensional flat 
case,  we obtain with Maple:
$$
B_{4,\flat}(f,h)
= -4\Bigl(\llucovd{f}{i}{j}{j}\ucovd{h}{i} + \llucovd{h}{i}{j}{j}\ucovd{f}{i}\Bigr)
   -4\llcovd{f}{i}{j} \uucovd{h}{i}{j} -2\lucovd{f}{i}{i} \lucovd{h}{j}{j},
$$
and in the 6 dimensional flat case:
\begin{align*}
&B_{6\,\flat}(df,dh) \\
&= 12\,(\lcovd{f}{i}\ululucovd{h}{i}{j}{j}{k}{k} + \lcovd{h}{i}\ululucovd{f}{i}{j}{j}{k}{k})
 + 24\,(\llcovd{f}{i}{j}\uulucovd{h}{i}{j}{k}{k} + \llcovd{h}{i}{j}\uulucovd{f}{i}{j}{k}{k}) \\
&\quad + 6\,(\lucovd{f}{i}{i}\lulucovd{h}{j}{j}{k}{k} + \lucovd{h}{i}{i}\lulucovd{f}{j}{j}{k}{k})
       + 24\,\llucovd{f}{i}{j}{j}\ulucovd{h}{i}{k}{k}
       + 16\lllcovd{f}{i}{j}{k}\uuucovd{h}{i}{j}{k},
\end{align*}
Here each summand is explicitly symmetric on $f$ and $h.$

\begin{remark}
\label{universalflat}
By looking at the construction of $B_n$ in the flat case and by Theorem~\ref{lemmaomeganflat}, 
we know that the coefficients of $B_{n\,\flat}$ are obtained from the Taylor expansion of the 
function $\psi(\xi,\xi)$ and integration on the sphere $||\xi||=1.$ Hence they are universal in the flat 
case.  Furthermore, since the volume of the sphere $\bbS^{n-1}$ is $2\pi^{n/2}/\Ga(n/2),$ any integral 
of a polynomial in the $\xi_i$'s with rational coefficients will be a rational multiple of $\pi^{n/2}.$ 
To avoid the factor $\pi^{n/2},$ we have assumed in \nnn{wresP} $d \,\xi$ to be normalized.
We conclude from Theorem~\ref{lemmaomeganflat} that the coefficients of $B_{n\,\flat}(f,h)$ are all 
rational in $n.$ 
\end{remark}

In the particular case in which $\widehat g = e^{2\eta}g$ with $g$ the flat metric on $M,$ we have the
conformal change equation for the Ricci tensor:
\begin{equation}
\label{invariantize} \llcovd{\eta}{i}{j} = - \Rho_{ij} - \lcovd{\eta}{i}\,\lcovd{\eta}{j} +
\frac{1}{2}\lcovd{\eta}{k}\,\ucovd{\eta}{k}\,g_{ij},
\end{equation}
which allows to replace all higher derivatives on $\eta$ with terms with the Ricci tensor.
In this way, using \nnn{invariantize} it is possible to obtain the
expression for $B_n$ in the conformally flat case out of its expression in the flat metric.  
We exemplified this process in \cite{mio2} for the six dimensional case.

\begin{remark}
The universality of the expression for $B_{n\,\flat}$ is preserved into the universality of 
the expression for $B_n$ in the conformally flat case, $B_{n\,\confflat}.$  In the next section we 
will see in general that there 
is a universal expression for $B_n.$  Furthermore, using \nnn{invariantize}
we conclude that also the coefficients of $B_{n\,\confflat}$ are rational in $n.$
\end{remark}

%%-----> A recursive approach to the symbol of $F$:-------------------------------------------------

\section{A recursive approach to the symbol of $F$}
\label{thesymbolofF}

In this section we give a recursive way of computing $\s_{-j}(F)$ in some given local charts, by
considering $F=2\cD-1$ with $\cD$ the orthogonal projection on the image of $d$ inside
$L^2(M,\Lambda^{n/2}_\mathbb{C} T^*M)$.  Both $F$ and $\cD$ are pseudo-differential operators of order
0 with $\s_0(F) = 2\s_0(\cD) - I$ and $\s_j(F) = 2 \s_j(\cD)$ for each $j < 0$. We shall use the
symbols of the differential operators of order 2, $d\d$ and $\Delta$.  The purpose of this section is to 
understand the relation in between $B_n$ and the curvature tensor.  

\begin{remark}
Even though the different 
$\s_{-j}(F)$ are not diffeomorphic invariant, an explicit understanding of their expressions in a given 
coordinate chart can be of optimal use, an example is present in \cite{Kastler} and \cite{KaWa} where 
using normal coordinates and recursive computation of symbols they were able to prove that the action 
functional, that is, the Wodzicki residue of $\Delta^{-n/2+1},$ is proportional to the integral of the 
scalar curvature by a constant depending on $n:$
$$
\Wres \Delta^{-n/2+1} = \frac{(n/2-1)\Omega_n}{6}\int_M \Sc |v_g|,
$$
with $\Omega_n$ the volume of the standard $n$-sphere, and $|v_g|$ the 1-density associated to 
the normalize volume form of $M.$
\end{remark}

From Lemma~{2.4.2}~\cite{Gilkey}, it is possible to understand the relation between the
$\s_j^\Delta$ and the metric at a given point $x.$  The same reasoning can be applied to
the operator $d\d-\d d$ and by addition the same is true for $d\d.$  We represent $\s(\Delta),$ 
$\s(d\d-\d d),$ 
and $\s(d\d)$ as $\s(\Delta) = \s_2^\Delta + \s_1^\Delta + \s_0^\Delta,$  
$\s(d\d-\d d) = \s_2^{d\d-\d d} + \s_1^{d\d-\d d} + \s_0^{d\d-\d d},$ and 
$\s(d\d) = \s_2^{d\d} + \s_1^{d\d} + \s_0^{d\d}$ with each $\s_j$ the component of homogeneity $j$ on the 
co-variable $\xi$.  The conclusion reads: $\s_2$ only invokes $g_{ij}(x),$
$\s_1$ is linear in the first partial derivatives of the metric at $x$ with coefficients depending smoothly
on the $g_{ij}(x)$, and $\s_0$ can be written as a term linear in the second partial
derivatives of the metric at $x$ with coefficients depending smoothly on the $g_{ij}(x)$ plus
a quadratic term linear in the first partial derivatives of the metric at $x.$  

Since in terms of a Riemannian metric compatible with
the conformal structure of $M$ we have $\Delta \cD=d\d$ we know $\s(\Delta \cD) = \s(d\d).$  Thus
the formula for the total symbol of the product of two pseudo-differential operators implies
\begin{align*}
\s_2^{d\d} + \s_1^{d\d} + \s_0^{d\d} &= \s(d\d) = \s(\Delta \cD)
\sim \sum\frac{1}{\a!}\del^\a_\xi \s(\Delta) D^\a_x (\s(\cD)) \\
&\sim \sum \frac{1}{\a !} \del^\a_\xi(\s_2^\Delta + \s_1^\Delta
+ \s_0^\Delta) D^\a_x(\s_0^\cD + \s_{-1}^\cD+ \s_{-2}^\cD + \cdots).
\end{align*}
Expanding the right hand side into sum of terms with the same homogeneity we conclude:

\begin{lemma}
\label{lemmasymbolofcD} 
In any given system of local charts, we can express the total symbol of
$\cD$, $\s(\cD) \sim \s_0^\cD + \s_{-1}^\cD + \cdots$ in a recursive way by the formulae:
\begin{align*}
\s_0^\cD &= (\s_2^\Delta)^{-1} \s_2^{d\d}, \\
\s_{-1}^\cD &= (\s_2^\Delta)^{-1}\Bigl(\s_1^{d\d} - \s_1^\Delta \s_0^\cD
 - \sum_{|\a|=1} \del^\a_\xi(\s_2^\Delta) D^\a_x(\s_0^\cD) \Bigr), \\
\s_{-2}^\cD &= (\s_2^\Delta)^{-1}\Bigl(\s_0^{d\d} - \s_1^\Delta \s_{-1}^\cD - \s_0^\Delta \s_0^\cD 
\\ &\quad
   - \sum_{|\a|=1}\left(\del^\a_\xi(\s_2^\Delta) D^\a_x(\s_{-1}^\cD) 
   + \del^\a_\xi(\s_1^\Delta) D^\a_x(\s_0^\cD)\right)
- \sum_{|\a|=2} \frac{1}{\a!} \del^\a_\xi(\s_2^\Delta)D^\a_x(\s_0^\cD) \Bigr),
\\
\s_{-r}^\cD &= -(\s_2^\Delta)^{-1}\Bigl( \s_1^\Delta \s_{-r+1}^\cD + \s_0^\Delta \s_{-r+2}^\cD  
+\sum_{|\a|=1}\del^\a_\xi(\s_2^\Delta) D^\a_x(\s_{-r+1}^\cD)
\\ &\qquad\qquad\quad
+\sum_{|\a|=1}\del^\a_\xi(\s_1^\Delta) D^\a_x(\s_{-r+2}^\cD)
+ \sum_{|\a|=2} \frac{1}{\a!} \del^\a_\xi(\s_2^\Delta)D^\a_x(\s_{-r+2}^\cD)\Bigr),
\end{align*}
for every $r \geq 3.$  Thus the total symbol of $F$, $\s^F \sim \s_0^F + \s_{-1}^F + \cdots,$ can
be recover from the relations $\s_0^F = 2\s^\cD_0 - I$ and $\s_{-k}^F = 2 \s_{-k}^\cD$ for $k\leq 1.$  
\end{lemma}

\begin{lemma}
\label{universalityofomega}
The bilinear differential functional $B_n$ on Theorem~\ref{theoremonOmegan} has a universal expression as a 
polynomial on $\Na$, $R,$ $df,$ and $dh$ with coefficients rational on $n.$
\end{lemma}

\begin{proof}
By choosing the coordinates to be normal coordinates, we
can assume the following: $g_{ij}(x)={\d_K}_{ij},$ that the partial derivatives of the metric vanish at $x,$
and that any higher order partial derivative of the metric at $x$ are expressed as polynomials
on $\Na$ and $R$ (see for example Corollary~2.9~\cite{Gray}).
By Lemma~\ref{lemmasymbolofcD} each $\s_{-k}^F(x,\xi)$ has a polynomial expression on $\Na$ and $R$ 
at the point $x.$  Every product
$\del^{\a}_\xi(\s^F_{-i}) \del^{\b}_\xi(D^{\ga}_x(\s^F_{-j}))$ will have a polynomial expression on 
$\Na$ and $R.$  These properties are preserved after integration of $\xi$ on $||\xi||=1$ and hence
$$
B_n(f,h) := 
\int_{||\xi||=1} 
\sum \frac{D^\b_x(f) D^{\a''+\d}_x(h)}{\a'!\a''!\b!\d!}
\tr\biggl(
\del^{\a'+\a''+\b}_\xi(\s^F_{-i}) \del^{\d}_\xi(D^{\a'}_x(\s^F_{-j})) \biggr)\,d \xi
$$
with sum taken over $|\a'|+|\a''|+|\b|+|\d|+i+j=n,$ $|\b|\geq 1,$ $|\d|\geq 1$ has an expression as a
polynomial in the ingredients $\Na$, $R,$ $df,$ and $dh.$  The explicit expression from 
Lemma~\ref{lemmasymbolofcD} for each $\s_{-k}^F$ in terms of $\s_i^\Delta$ and $\s_j^{d\d}$ depends on 
the metric, the dimension, and the local coordinates since the different $\s_j$ are not diffeomorphic 
invariant, but \nnn{wresP} and hence $B_n$ is independent of local representations.  We conclude 
that the expression we obtain for $B_n$ is universal.

This independence of local representations allows as to conclude the rationality of the coefficients 
of $B_n$ as follows. By Lemma~\ref{lemmasymbolofcD} each $\s_{-j}^\cD$ is a polynomial expression 
in the $\del_x\s_k^\Delta,$ $\del_x\del_\xi\s_k^\Delta,$ and  $\del_x\s_k^{d \d},$ with 
coefficients rational in $n.$  
By choosing the local coordinates to be normal coordinates, we can assure that all the expressions of 
these symbols as polynomials on $\Na$ and $R$ are of coefficients rational in $n.$    
\end{proof}

%%-----> A filtration by degree:-------------------------------------------------------------------

\subsection{A filtration by degree}
\label{afiltrationbydegree}

If we fix $h$ (or $f$), $B_n(f,h)$ is a differential operator or order $n-1$ on $f$ (or $h$).  It is
acting on smooth functions and producing smooth functions.  Since $\widehat{B_n} = e^{-n\eta} B_n$ 
when we transform the metric conformally by $\widehat g =e^{2\eta} g,$ we have
that $B_n$ is of level $n,$ (see \cite{Bra4}) that is, the expression $B_n(f,h)$ is a sum of 
homogeneous polynomials
in the ingredients $\Na^\a f,$ $\Na^\b h,$ and $\Na^\ga R$ for multi-indices $\a,\b,$ and $\ga,$
in the following sense, each monomial must satisfies the homogeneity condition given by the
rule: 
$$ 
\mbox{twice~the~appearances~of~$R$~} + \mbox{~number~of~covariant~derivatives} = n 
$$
where for covariant derivatives we count all of the derivatives on $R,$ $f,$ and $h,$ and any
occurrence of $W,$ $Rc,$ $\Rho,$ $Sc,$ or $J$ is counted as an occurrence of~$R.$ Furthermore, by the
restriction $|\b|\geq 1$ and $|\d|\geq 1$ in Lemma~\ref{lemmas_nSS}, we know that $B_n(f,h)$ is made of
the ingredients $\Na^\a df,$ $\Na^\b dh,$ and $\Na^\ga R.$

By closing under addition, we denote by $\cP_n$ the space of these polynomials.
For a homogeneous polynomial $p$ in $\cP_n,$ we denote by $\deg_R$ its degree in~$R$ and by $\deg_\Na$
its degree in~$\Na.$  In this way, $2\deg_R + \deg_\Na = n,$ with $\deg_\Na \geq 2$ and hence 
$2\deg_R \leq n-2$ for $B_n(f,h).$  We say that $p$ is in $\cP_{n,r}$ if $p$ can be written as a sum 
of monomials with $\deg_R\geq r,$ or equivalently, $\deg_\Na \leq n - 2r.$ We have a filtration by degree:
$$
\cP_n = \cP_{n,0} \supseteq \cP_{n,1} \supseteq \cP_{n,2}
\supseteq\cdots\supseteq \cP_{n,\tfrac{n-2}{2}},
$$
and $\cP_{n,r} = 0$ for $r> (n-2)/2.$ There is an important observation to make. An expression
which a priori appears to be in, say $\cP_{6,1},$ may actually be in a subspace of it, like $\cP_{6,2}.$
For example,
\begin{equation}
\label{forfiltration}
\underbrace{\lcovd{f}{i} \lllcovd{h}{j}{k}{l} W^{ijkl}}_{\in\,\cP_{6,1}} =
\underbrace{\lcovd{f}{i}\lcovd{h}{j} \Rho_{kl} W^{ikjl}
+ \frac{1}{2}\lcovd{f}{i}\lcovd{h}{j} W^i{}_{klm}W^{jklm}}_ {\in \,\cP_{6,2}},
\end{equation}
by reordering covariant derivatives and making use of the symmetries of the Weyl tensor. 

In the particular case $n=4,$ $k_R$ can be 0 or~1, hence $B_4$ can be written as
$$
B_{4}(f,h) = \underbrace{p_4(df,dh)}_{\in\,\cP_{4,0}} + \underbrace{p_{R,4}(df,dh)}_{\in\,\cP_{4,1}}
$$
where $p_{R,4}(df,dh)$ is a trilinear form on $R,$ $df,$ and $dh.$  Explicitly,
$$
p_4(df,dh) = B_{4\,\flat}(f,h)
\hbox{~~and~~}
p_{R,4}(df,dh)= 8 \lcovd{f}{i} \ucovd{h}{i} J.
$$
In the 6-dimensional case, $k_R \in \{0,1,2\}$ thus
$$
B_{6}(f,h) =
\underbrace{p_6(df,dh)}_{\in \cP_{6,0} / \cP_{6,1}}
+ \underbrace{p_{6,1}(df,dh)}_{\in \cP_{6,1} / \cP_{6,2}}
+ \underbrace{p_{6,2}(df,dh)}_{\in\cP_{6,2}},
$$
with
$$
p_{6,1} = p_{6,R}+ p_{6,R'} + p_{6,R''},
$$
\begin{itemize}
\item[-]
$p_6(df,dh) = B_{6\,\flat}(f,h),$
\item[-]
$p_{6,R}(df, dh),$ polynomial on $\Na^a df,$ $\Na^b dh,$ and $R,$
\item[-]
$p_{6,R'}(df, dh),$ polynomial on $\Na^a df,$ $\Na^b dh,$ and $\Na R,$
\item[-]
$p_{6,R''}(df, dh),$ polynomial on $\Na^a df,$ $\Na^b dh,$ and $\Na \Na R,$ and
\item[-]
$p_{6,2}(df, dh),$ polynomial on $\Na^a df,$ $\Na^b dh,$ and $R R.$
\end{itemize}

From the previous expressions, it is evident that there exists a sub-filtration inside each
$\cP_{n,l}$ for $l\geq 1.$ Such a filtration is more complicated to describe in higher
dimension because of the presence of terms like $\Na^a R \Na^c R\cdots.$  Also it is important to 
note that
$p_n(df,dh)$ is precisely $B_{n,\flat}(f,h),$ i.e. the flat version of $B_n$ coincide with the expression
in the filtration without curvature terms.  

%%-----> $P_n$ and the Wodzicki residue:---------------------------------------------------------------

\section{$P_n$ and the Wodzicki residue}
\label{PnandWres}

Now we are ready to prove our main result:
\begin{proof}{\bf Of Theorem~\ref{theoremonP_n}}.
Since $B_n(f,h)$ is a differential functional on $f$ and $h$, Stokes' theorem applied to
$\int_M B_n(f,h) \,d x$ leads to the expression $\int_M f P_n(h) \,d x,$ where $P_n$ is a
differential operator on $h.$  Uniqueness of $P_n$ follows from the arbitrariness of $f.$

Formally selfadjointness is a consequence of the symmetry of $B_n,$ that is:
$$
\int_M f P_n(h) \,d x = \int_M B_n(f,h) \,d x
                        = \int_M B_n(h,f) \,d x
                        =\int_M P_n(f) h \,d x.
$$

Note that
\begin{align*}
\int_M f e^{-n\eta} P_n(h) \,d x
&= \int_M B_n(fe^{-n\eta},h) \,d x = \int_M \widehat{B_n}(fe^{-n\eta},h) \,\widehat{d x} \\
&= \int_M f e^{-n\eta} \widehat{P_n}(h) \,\widehat{d x} =\int_M f \widehat P_n(h) \,d x
\end{align*}
for every $f \in C^\infty(M).$ (ii) follows from the arbitrariness of $f.$

From Lemma~\ref{universalityofomega} there is a universal expression for $B_n$ and
Stoke's theorem does not affect this universality nor the rationality of its coefficients, thus there is 
a universal expression for $P_n$ as a polynomial on $\Na$ and $R,$
with coefficients rational on $n,$ and (iii) is proved.

Next we prove (iv).  Because of the filtration, we can write  $B_n(f,h)$ in the form 
$B_n(f,h) = p_n(df,dh) + p_{n,R}(df,dh)$ where any possible curvature term is in the $p_{n,R}$ part 
and those terms in $p_n$ are of the form 
\begin{equation}
\label{fhRfirstkind}
f_{;\,j_1}{}^{j_1}\cdots{}_{j_s}{}^{j^s}{}_{i_1 \cdots i_r} 
\,
h_{;\,k_1}{}^{k_1}\cdots{}_{k_t}{}^{k_t i_1 \cdots i_r} 
\end{equation}
with $2r+2s+2t=n.$ Each time we apply Stokes' theorem
to such a term and reorder the indices, we obtain an expression like:
$$ 
(\pm)\int f \Delta^{n/2}(h)\,d x + \int f \lot(h,R)\, d x 
$$ 
where $\lot(h,R)$ represents a differential operator on $h$ which is polynomial on $\Na$ and $R$ 
with order on $h$ lower than $n.$ It follows that $P_n(h) = c_n \Delta^{n/2}(h) + \lot$ with $c_n$ a 
universal constant and $\lot$ a sum of lower order terms. 

To prove (v) we take a closer look at $B_n(f,h).$  Each time we apply Stokes' theorem
to a term of the form \nnn{fhRfirstkind} for which $r\geq 1,$ 
we can reorder the indices in such a way that we have an expression like:
$$
(\pm)\int f \Delta^{n/2}(h)\,d x + \int f \bigl(\lot(h_i,R)\bigr)_{;\,}{}^i\, d x
$$
where now $\lot(h_i,R)$ represents a sum of contractions of products of the form 
$\Na^\a h\times$ (some curvature term), for some multi-index $\a$ with $i$ present in $\a.$  
Each time we apply Stokes' theorem to a term of the form \nnn{fhRfirstkind} for which $r=0,$  
we can reorder the indices in such a way that we have an expression like:
\begin{equation}
\label{fhRthirdkind}
(\pm)\int f \Delta^{n/2}(h)\,d x + \int f \bigl(\lot(h,R)\bigr)_{;\,i}{}^i\, d x
\end{equation}
where $\lot(h,R)$ represents a sum of lower order terms.
Each time we apply Stokes' theorem to any term of $B_n(f,h),$ not of the form 
\nnn{fhRfirstkind}, that is to say with some curvature term, we can reorder the indices in such a 
way that we have an expression of the form:
\begin{equation}
\label{fhRfourthkind}
\int f\bigl(\lot(h_i,R)\bigr)_{;}{}^i\,d x 
+ \int f\bigl(\lot(h,R_i)\bigr)_{;}{}^i\,d x
\end{equation}
where $\lot(h,R_i)$ represents a sum of contractions of products of the form
$\Na^\a h\times$ (some curvature term), for some multi-index $\a$ in such a way that $i$ is present in 
one of the curvature terms (not in $\a$).

Applying the Leibniz rule to the first appearance of $i$ in the second summand of \nnn{fhRthirdkind} 
produces more terms of the same form as \nnn{fhRfourthkind}.  We conclude: 
$$
P_n(h) = c_n\Delta^{n/2}(h) 
+ \sum \bigl(\lot(h_{i},R)\bigr)_{;\,}{}^{i} 
+ \sum \bigl(\lot(h,R_{i})\bigr)_{;\,}{}^{i}.
$$

Next, consider the operator $S'_n$ acting on exact 1-forms:
$$
S'_n :dh = h_{;\,i}\,d x^i \mapsto - \Bigl(\sum \lot(h_{i},R) + \sum \lot(h,R_{i})\Bigr)\,d x^i.
$$
It follows that $P_n(h) = c_n\Delta^{n/2}h + \delta S'_n d(h).$
Hence $P_n = \delta S_n d$ where $S_n$ is a constant multiple of $\Delta^{n/2-1} + \lot$ or 
$(d \delta)^{n/2-1} + \lot,$ where any curvature term is absorbed by the $\lot$ part.

For (vi), we will prove that for any $k\in C^\infty(M)$ we have
\begin{equation}
\label{toprove} 
\int_M k \bigl(P_n(fh) - fP_n(h) - hP_n(f)\bigr) \,d x = -2\int_M k\,B_n(f,h)\,d x
\end{equation}
and the result will follow from the arbitrariness of $k.$  By \nnn{omegatoP} the left hand side of
\nnn{toprove} is equal to:
\begin{align*}
& \int_M B_n(k,fh) - B_n(kf,h) - B_n(kh,f) \,d x \\
&= \Wres\bigl([F,k][F,fh]-[F,kf][F,h]-[F,kh][F,f]\bigr) \\
&= \Wres\bigl(FkFfh - FkfhF - kFFfh + kFfhF \\
&\qquad\quad -FkfFh + FkfhF + kfFFh - kfFhF \\
&\qquad\quad -FhkFf + FhkfF + hkFFf - hkFfF \bigr) \\ 
&= -2\Wres(kFfFh - kFfhF - kfFFh + kfFhF) \\
&= -2\Wres\bigl(k[F,f][F,h]\bigr) = -2 \int_M k\, B_n(f,h)\,d x,
\end{align*}
where we have used the commutativity of the algebra $C^\infty(M),$ the
linearity and the trace property of $\Wres,$ and the property $F^2=1.$
\end{proof}

\begin{remark}
The value of $c_n$ in Theorem~\ref{theoremonP_n}~(iii) gets determined in the flat case since no terms 
with curvature affect it.  In the 4 dimensional case $c_4=2$ and in the 6 dimensional case $c_6 =-4.$
\end{remark}

The relation between the {\sf GJMS} operators in the general even dimensional case and the $P_n$ 
operators described in this work still unknown to the author at this time.  
Nevertheless, in the flat case, both operators coincide up to a constant multiple with a power of 
the Laplacian, $\Delta^{n/2},$ and both operator 
enjoy the same conformal property 
$$
\widehat P=e^{-n\eta} P
$$
for $\widehat g = e^{2\eta} g.$ Hence, they must also coincide, 
up to a constant multiple, in the conformally flat case.
\begin{proposition}
\label{lastprop} 
In the even dimensional case, inside the conformally flat class of metrics, the critical {\sf GJMS} 
operator and the operator $P_n$ described in this work coincide up to a constant multiple.
\end{proposition}

%%-----> REFERENCES: ----------------------------------------------------------------------------

\obeylines{
Department of Mathematics, Purdue University.  ugalde@math.purdue.edu
}

\end{document}